\documentclass{elsarticle}

\usepackage{lineno,hyperref}
\usepackage{amsmath}
\usepackage{amsfonts,amssymb}
\usepackage{esint}
\usepackage{subfigure}
\usepackage{wasysym}
\usepackage{amsthm}
\usepackage{multirow}
\usepackage{algorithm}
\usepackage{algorithmic}

\newtheorem{remark}{Remark}

\modulolinenumbers[5]

\journal{Elsevier}

\begin{document}

\begin{frontmatter}

\title{An alternative reconstruction for WENO schemes with adaptive order}


\author[add1]{Hua Shen\corref{mycorrespondingauthor}}
\cortext[mycorrespondingauthor]{Corresponding author}
\ead{huashen@uestc.edu.cn}
\address[add1]{School of Mathematical Science, University of Electronic Science and Technology of China, Chengdu, Sichuan 611731, China}

\begin{abstract}
We propose an alternative reconstruction for weighted essentially non-oscillatory schemes with adaptive order (WENO-AO)
for solving hyperbolic conservation laws.
The alternative reconstruction has a more concise form than the original WENO-AO reconstruction.
Moreover, it is a strictly convex combination of polynomials with unequal degrees.
Numerical examples show that the alternative reconstruction maintains the accuracy and robustness of the WENO-AO schemes.
\end{abstract}


\end{frontmatter}


\section{Introduction}\label{sec:intro}
The solution to a hyperbolic conservation law commonly contains both discontinuities and smooth regions with sophisticated structures,
so numerical schemes with high-order convergence and excellent shock-capturing properties are desired for solving general hyperbolic conservation laws.
The essentially non-oscillatory (ENO) schemes and weighted ENO (WENO) schemes are state-of-the-art such schemes
and become mainstream in many applications.

The history of ENO schemes can date back to 1980s.
Harten and his coworkers \cite{Harten1987ENO} constructed uniformly high-order accurate ENO schemes by slightly relaxing the TVD constraint \cite{Harten1983TVD}.
ENO schemes achieve uniformly high-order accuracy and the non-oscillatory feature by the strategy of choosing the smoothest stencil from several local candidates.
Shu and Osher \cite{Shu1988EfficientENO, Shu1989EfficientENOII} significantly improved the efficiency of ENO schemes
by using TVD Runge-Kutta time discretizations and numerical fluxes reconstructions instead of cell average reconstructions.
Based on ENO schemes, Liu \emph{et al}. \cite{Liu1994WENO} proposed WENO schemes by assigning proper weights to the local candidate stencils.
Jiang and Shu \cite{Jiang_Shu1996WENO} provided a general framework for designing smooth indicators and weights for WENO schemes
to achieve the optimal order in smooth regions while maintaining the non-oscillatory feature at discontinuities.
Borges \emph{et al}. \cite{Borges2008WENO_Z} further improved the accuracy of WENO schemes by redesigning the nonlinear weights.
Besides, many extended ENO/WENO schemes are proposed and a lot of shock capturing high-order schemes were constructed in the light of ENO/WENO ideas,
such as the monotonicity preserving WENO schemes \cite{Balsara2000Monotonicity_WENO}, the weighted compact nonlinear schemes \cite{Deng2000WCNS},
the compact central WENO schemes \cite{Bianco1999CENO, Levy1999CWENO, Levy2000CompactCWENO}, the Hermite WENO schemes \cite{Qiu2004HermiteWENO, Qiu2005HermiteWENO2D}, the  WENO-ADER schemes \cite{Toro2002ADER, Titarev2004WENO_FV, Balsara2013ADER}, the $P_NP_M$ schemes \cite{Dumbser2008PNPM, Dumbser2009PNPM, Dumbser2010PNPM}, and so forth.

Recently, the WENO schemes with adaptive order (WENO-AO) attracted a lot of interest.
WENO-AO schemes started from the fifth-order WENO-ZQ scheme \cite{Zhu2016NewWENOFD}, and was further developed by Balsara \emph{et al.} \cite{Balsara2016WENO_AO}.
The building block of WENO-AO schemes is the combination of polynomials with different degrees that was originally proposed by Levy \emph{et al}. \cite{Levy2000CompactCWENO}.
The main advantage of WENO-AO schemes is that they can adjust their order adaptively depending on the smoothness of the solutions
thanks to the building block.
Therefore, it is easy to construct numerical schemes with multi-resolution \cite{Balsara2016WENO_AO, Balsara2020WENO_AO_Unstructure, Zhu2018NewMultiResolutionWENO, Zhu2019NewMultiResolutionWENO_Tri, Zhu2020NewMultiResolutionWENO_Tet}.
In the original reconstruction \cite{Levy2000CompactCWENO, Zhu2016NewWENOFD, Balsara2016WENO_AO},
the high-order polynomial needs to subtract the low-order polynomials in order to recover the optimal high-order polynomial in smooth regions.
In doing so, the reconstructed polynomial is not a strictly convex combination of the low-order and high-order polynomials.
Balsara \emph{et al.} \cite{Balsara2016WENO_AO}  already realized this issue when constructing the WENO-AO(7,5,3) scheme.
Although this slight imperfection seems not to affect the practical behavior,
a purist might still seek for a cure.
In this note, we propose an alternative reconstruction which ensures the convexity property
and maintains the accuracy and efficiency of WENO-AO schemes.

\section{A brief review of finite difference WENO-AO schemes}\label{SEC:Original_WENO_AO}
We consider the one-dimensional scalar conservation law,
\begin{equation}\label{Eq:1DHCL_Eq}
    \frac{\partial u}{\partial t}+\frac{\partial f(u)}{\partial x}=0, x\in[x_L,x_R],t\in[0,\infty).
\end{equation}

The spatial domain $[x_L,x_R]$ is discretized into uniform intervals
by $x_i=x_L+(i-1)\Delta x$ ($i=1 \text{ to } N+1$), where $\Delta x=(x_R-x_L)/N$.
The finite difference WENO schemes can be expressed in the following conservative form
\begin{equation}\label{Eq:1DWENOFD}
  \frac{du_i(t)}{dt}=\mathcal{L}(u_i(t))=-\frac{\hat{f}_{i+1/2}-\hat{f}_{i-1/2}}{\Delta x},
\end{equation}
where the numerical flux $\hat{f}_{i\pm1/2}$ approximates the function $h(x)$,
that is implicitly defined by
$f(u(x))=\frac{1}{\Delta x}\int_{x-\Delta x/2}^{x+\Delta x/2}h(\xi)d\xi$ \cite{Shu1988EfficientENO}, at $x_{i\pm1/2}$.

In order to enhance the robustness of the scheme,
we usually introduce the upwind mechanism by choosing an upwind-biased stencil according to the direction of the wave propagation,
i.e., the sign of $\frac{df(u)}{du}$.
For a general flux, the sign of $\frac{df(u)}{du}$ is not constant,
but it can be split into two parts as
\begin{equation}\label{Eq:Flux_Split}
  f(u)=f^+(u)+f^-(u),
\end{equation}
where $\frac{df^+(u)}{du}\ge0$ and $\frac{df^-(u)}{du}\le0$.
For example, the global Lax–Friedrichs flux splitting method splits the flux as
\begin{equation}\label{Eq:LF_Flux_Split}
  f^\pm(u)=\frac{1}{2}(f(u)\pm\alpha u),
\end{equation}
where $\alpha=\mbox{max}\left|\frac{df(u)}{du}\right|$
and the maximum is taken over the whole computational domain.

After the flux splitting, we respectively construct polynomials to approximate $\hat{f}^\pm(u)$ on several sub-stencils of an upwind-biased large stencil,
and then use a nonlinear function of the stencils' smoothness indicators to combine the constructed polynomials.
Balsara \emph{et al.} \cite{Balsara2016WENO_AO} provided an efficient way
to construct a polynomial and to calculate the corresponding smoothness indicator
on a given stencil by using Legendre basis.
In classical WENO reconstructions \cite{Jiang_Shu1996WENO, Borges2008WENO_Z},
the sub-stencils have an equal size, and so do the associated polynomials.
By carefully design the weights for the sub-stencils,
the reconstructed polynomial recovers the optimal high-order polynomial on the large stencil in smooth regions
and tends to the low-order polynomial on the smoothest sub-stencil.
However, the reconstructed polynomial cannot always recover the optimal high-order polynomial,
especially for multidimensional reconstructions.
The WENO-AO reconstruction fixes this issue via a non-linear hybridization between the optimal high-order polynomials and low-order polynomials.
For example, the WENO-AO(5,3) reconstruction is expressed as \cite{Balsara2016WENO_AO}
\begin{equation}\label{Eq:WENO_AO_5_3}
\begin{split}
   P^{\mbox{AO(5,3)}}(x)&=\frac{\bar{w}_3^{r5}}{\gamma_3^{r5}}\left(P_3^{r5}(x)-\gamma_1^{r3}P_1^{r3}(x)-\gamma_2^{r3}P_2^{r3}(x)-\gamma_3^{r3}P_3^{r3}(x)\right) \\
     & +\bar{w}_1^{r3}P_1^{r3}(x)+\bar{w}_2^{r3}P_2^{r3}(x)+\bar{w}_3^{r3}P_3^{r3}(x),
\end{split}
\end{equation}
where $P_3^{r5}(x)$ is the optimal fourth-order polynomial, $P_k^{r3}(x)$ $(k=1,2,3)$ are second-order polynomials,
and $\gamma$ and $\bar{w}$ are corresponding linear and nonlinear weights.
The specific forms of the polynomials and the weights can be found in \cite{Balsara2016WENO_AO}.
It is easy to verify that the reconstructed polynomial $P^{\mbox{AO(5,3)}}(x)$ recovers the optimal polynomial $P_3^{r5}(x)$ when $\bar{w}_n=\gamma_n$.
Using the above idea, we can construct numerical schemes with multi-resolution \cite{Balsara2016WENO_AO, Balsara2020WENO_AO_Unstructure, Zhu2018NewMultiResolutionWENO, Zhu2019NewMultiResolutionWENO_Tri, Zhu2020NewMultiResolutionWENO_Tet}.

It is straightforward to extend the finite difference WENO schemes to multi-dimensional cases,
because the same reconstruction procedure for approximating the numerical fluxes
can be implemented in a dimension-by-dimension manner.
In order to eliminate spurious oscillations as much as possible,
the reconstruction procedure is applied in the characteristic space when solving hyperbolic systems.

As for the time integration,
we adopt the method of lines which first approximates the numerical fluxes by WENO reconstructions
and uses the third-order TVD Runge-Kutta method \cite{Shu1988EfficientENO} to solve the system of ordinary differential equations,
i.e., $\frac{du_i(t)}{dt}=\mathcal{L}(u_i(t))$.
\section{An alternative WENO-AO reconstruction}\label{SEC:Alternative_WENO_AO}
As we can see, the first part on the right hand side of Eq. (\ref{Eq:WENO_AO_5_3})
contains negative coefficients for $P_k^{r3}(x)$ $(k=1,2,3)$.
Therefore, $P^{\mbox{AO(5,3)}}(x)$ is not a strictly convex combination of $P_3^{r5}(x)$ and $P_k^{r3}(x)$ $(k=1,2,3)$.
In order to fix this defect, we propose an alternative WENO-AO (WENO-AOA) reconstruction.
Without loss of generality, the WENO-AOA(5,3) reconstruction can be expressed as
\begin{equation}\label{Eq:WENO_AOA_5_3}
   P^{\mbox{AOA(5,3)}}(x)=\bar{w}_3^{r5}P_3^{r5}(x)+\bar{w}_1^{r3}P_1^{r3}(x)+\bar{w}_2^{r3}P_2^{r3}(x)+\bar{w}_3^{r3}P_3^{r3}(x),
\end{equation}
where $P_3^{r5}(x)$ and $P_k^{r3}(x)$ $(k=1,2,3)$ are the same as the original WENO-AO reconstruction.
The un-normalized weights are slightly modified from Borges \emph{et al.} \cite{Borges2008WENO_Z}
\begin{equation}\label{Eq:UnNormalized_Weights}
\begin{split}
   w_3^{r5}=\gamma_3^{r5}\left[1+\left(\frac{\tau}{\beta_3^{r5}+\epsilon}\right)^2\right],\quad  w_1^{r3}=\gamma_1^{r3}\left(\frac{\tau}{\beta_1^{r3}+\epsilon}\right)^2,\\
     w_2^{r3}=\gamma_2^{r3}\left(\frac{\tau}{\beta_2^{r3}+\epsilon}\right)^2, \quad w_3^{r3}=\gamma_3^{r3}\left(\frac{\tau}{\beta_3^{r3}+\epsilon}\right)^2,
\end{split}
\end{equation}
where $\epsilon=10^{-40}$ is used to avoid singularity,
$\gamma$ represents constant linear weights,
the smoothness indicator measuring the regularity of a $n$th order polynomial $p_n(x)$ on $[x_{i-1/2},x_{i+1/2}]$ is given by \cite{Jiang_Shu1996WENO}
\begin{equation}\label{Eq:Smoothness_Indicator}
  \beta_n=\sum_{l=1}^{n}\Delta x^{2l-1}\int_{x_{i-1/2}}^{x_{i+1/2}}\left(\frac{d^lp_n(x)}{dx^l}\right)^2dx,
\end{equation}
and the parameter $\tau$ is defined as
\begin{equation}\label{Eq:Tau}
  \tau=\frac{1}{3}\left(|\beta_3^{r5}-\beta_1^{r3}|+|\beta_3^{r5}-\beta_2^{r3}|+|\beta_3^{r5}-\beta_3^{r3}|\right).
\end{equation}
Finally, the normalized weights are given by
\begin{equation}\label{Eq:Normalized_Weights}
\begin{split}
   &\bar{w}_3^{r5}=\frac{w_3^{r5}}{w_3^{r5}+w_1^{r3}+w_2^{r3}+w_3^{r3}},\\
     &\bar{w}_k^{r3}=\frac{w_k^{r3}}{w_3^{r5}+w_1^{r3}+w_2^{r3}+w_3^{r3}}, \quad (k=1,2,3).
\end{split}
\end{equation}

The corresponding explicit expressions of the smoothness indicators can be found in \cite{Borges2008WENO_Z, Zhu2016NewWENOFD, Balsara2016WENO_AO}
and their Taylor series expansions at $x_i$ are given as
\begin{equation}\label{Eq:TaylorExpan}
\begin{split}
     & \beta_3^{r5}=f_i^{'2}\Delta x^2+\frac{13}{12}f_i^{''2}\Delta x^4+O(\Delta x^6), \\
     & \beta_1^{r3}=f_i^{'2}\Delta x^2+\left(\frac{13}{12}f_i^{''2}-\frac{2}{3}f'_if'''_i\right)\Delta x^4
     -\left(\frac{13}{6}f''_if'''_i-\frac{1}{2}f'_if''''_i\right)\Delta x^5+O(\Delta x^6), \\
     & \beta_2^{r3}=f_i^{'2}\Delta x^2+\left(\frac{13}{12}f_i^{''2}+\frac{1}{3}f'_if'''_i\right)\Delta x^4+O(\Delta x^6), \\
     & \beta_3^{r3}=f_i^{'2}\Delta x^2+\left(\frac{13}{12}f_i^{''2}-\frac{2}{3}f'_if'''_i\right)\Delta x^4
     +\left(\frac{13}{6}f''_if'''_i-\frac{1}{2}f'_if''''_i\right)\Delta x^5+O(\Delta x^6).
\end{split}
\end{equation}
For smooth solutions, it is easy to verify that
\begin{equation}\label{Eq:1DWmError}
    \bar{w}_k^{r3}\propto\begin{cases}
                    O\left(\Delta x^4\right), & \mbox{if } f'_i\ne 0, f'''_i\ne 0\\
                    O\left(\Delta x^2\right), & \mbox{if } f'_i=0, f''_i\ne 0,f'''_i\ne 0.
                  \end{cases}\quad (k=1,2,3).
\end{equation}
Since $\bar{w}_3^{r5}+\sum_{k=1}^{3}\bar{w}_k^{r3}=1$, the difference between the limited polynomial and the optimal high-order polynomial can be written as
   \begin{equation}\label{Eq:1DPolynomialError}
    \varepsilon=P^{\mbox{AOA(5,3)}}(x)-P_3^{r5}(x)=\sum_{k=1}^{3}\bar{w}_k^{r3}\left(P_k^{r3}(x)-P_3^{r5}(x)\right).
   \end{equation}
Because $P_3^{r5}(x)$ is a fifth-order approximation and $P_k^{r3}(x)$ $k=1,2,3$ are third-order approximations,
we have $P_k^{r3}(x)-P_3^{r5}(x)\propto\Delta x^3$.
Combining Eqs. (\ref{Eq:1DWmError}) and (\ref{Eq:1DPolynomialError}), we obtain
\begin{equation}\label{Eq:PolynomialError}
    \varepsilon\propto\begin{cases}
                    O\left(\Delta x^7\right), & \mbox{if } f'_i\ne 0, f'''_i\ne 0\\
                    O\left(\Delta x^5\right), & \mbox{if } f'_i=0, f''_i\ne 0,f'''_i\ne 0.
                  \end{cases}\quad (k=1,2,3).
\end{equation}
That is to say, $P^{\mbox{AOA(5,3)}}(x)$ does not violate the order of the optimal high-order polynomial $P_3^{r5}(x)$
in smooth regions including at local first-order extrema.

When the large stencil contains discontinuities,
we have $\tau\propto O\left(1\right)$, and $\beta_3^{r5}\propto O\left(1\right)$.
However, one of the low-order polynomial $P_k^{r3}(x)$ ($k$=1,2,3) is still smooth,
and we have $\beta_k^{r3}\propto O\left(\Delta x^2\right)$ ($k$=1 or 2 or 3).
Therefore, $\bar{w}_3^{r5}$ becomes very small, and $P^{\mbox{AOA(5,3)}}(x)$ approaches to the smooth low-order polynomial.
This property makes the reconstruction non-oscillatory at discontinuities.

\begin{remark}
  Comparing Eq. (\ref{Eq:WENO_AOA_5_3}) with Eq. (\ref{Eq:WENO_AO_5_3}),
  the proposed WENO-AOA reconstruction has a more concise form than the original WENO-AO reconstruction.
  More importantly, WENO-AOA reconstruction provided a strictly convex combination of polynomials with unequal degrees.
\end{remark}

\section{Numerical examples}
As suggested by Balsara \emph{et al}. \cite{Balsara2016WENO_AO},
we set the linear weights $\gamma_3^{r5}=\gamma_{Hi}$, $\gamma_1^{r3}=\gamma_3^{r3}=0.5(1-\gamma_{Hi})(1-\gamma_{Lo})$,
$\gamma_2^{r3}=(1-\gamma_{Hi})\gamma_{Lo}$, where $\gamma_{Hi}=\gamma_{Lo}=0.85$.
Without further illustration, we set CFL=0.5 for all cases.
\subsection{Linear advection}

\begin{table}
  \centering
  \begin{tabular}{|c|c|c|c|c|c|c|c|c|}
    \hline
      & \multicolumn{4}{|c|}{WENO-AO(5,3)} & \multicolumn{4}{|c|}{WENO-AOA(5,3)} \\
     \hline
      Mesh size       & $L_1$  &Order & $L_\infty$  &Order    & $L_1$  &Order & $L_\infty$  &Order    \\
     \hline
     1/25              &2.11E-6 &-     &3.38E-6      &-        &2.57E-6 &-     &4.11E-6      &- \\

     1/50             &6.67E-8 &4.98  &1.06E-7      &4.99     &8.11E-8 &4.99  &1.29E-7      &4.99 \\

     1/100             &2.09E-9 &5.00  &3.31E-9      &5.00     &2.54E-9 &5.00  &4.02E-9      &5.00 \\

     1/200             &6.66E-11 &4.97  &1.05E-10      &4.98     &8.02E-11 &4.99  &1.26E-10      &5.00 \\
     \hline

  \end{tabular}
  \caption{Numerical errors of the advection of the sinusoidal wave at $t=2$ computed by WENO-AO(5,3) and WENO-AOA(5,3) with different mesh sizes.}
  \label{Table:1DAdvectionConvergence}
\end{table}

We solve the scalar linear advection equation $\frac{\partial u}{\partial t}+\frac{\partial u}{\partial x}=0$ on the computational domain [-1, 1] with periodic boundaries.
The first case is the advection of the sinusoidal wave $u(x,0)=sin(\pi x)$.
Since the order of the time discretization does not match with the order of the space discretization,
we set $\Delta t=\Delta x^{5/3}$ to test the convergence rate of the numerical solutions.
Table \ref{Table:1DAdvectionConvergence} shows the average error $L_1$ and the maximum error $L_\infty$ at $t=2$
on different meshes induced by WENO-AO(5,3) and WENO-AOA(5,3).
We observe that the numerical errors of WENO-AO(5,3) and WENO-AOA(5,3) are quite close,
and both schemes achieve the optimal fifth-order convergence rate.

\begin{figure}
  \centering
  \subfigure[WENO-AO(5,3)]{
  \label{FIG:Linear_Advection_WENO_AO_5_3}
  \includegraphics[width=7 cm]{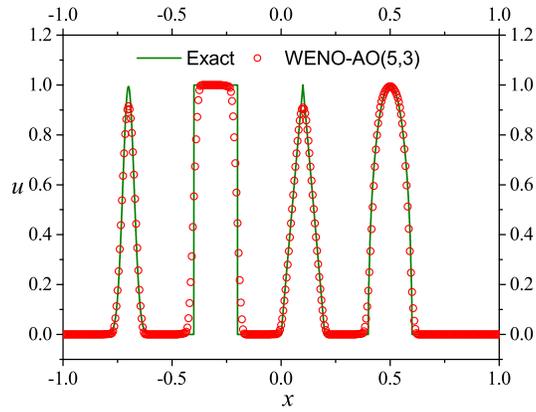}}
  \subfigure[WENO-AOA(5,3)]{
  \label{FIG:Linear_Advection_WENO_AOA_5_3}
  \includegraphics[width=7 cm]{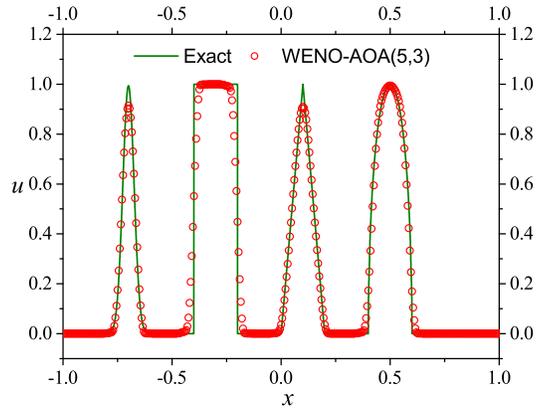}}
  \caption{The advection of a combination of Gaussians, a square wave, a sharp triangle wave,
  and a half ellipse arranged from left to right at $t=20$ calculated by WENO-AO(5,3) and WENO-AOA(5,3) with 401 mesh points.}
\label{FIG:Linear_Advection}
\end{figure}

The second case is the advection of a combination of Gaussians, a square wave, a sharp
triangle wave, and a half ellipse arranged from left to right.
This case was first proposed by Jiang and Shu \cite{Jiang_Shu1996WENO},
and the specific settings can be found therein.
The simulations are performed by using 401 mesh points and terminate at $t=20$.
Fig. \ref{FIG:Linear_Advection} shows the profiles of $u$ calculated by WENO-AO(5,3) and WENO-AOA(5,3).
We observe that the two schemes perform similarly for all kinds of solutions.

\section{Double Mach reflection problem}

\begin{figure}
  \centering
  \subfigure[WENO-AO(5,3)]{
  \label{FIG:DMR_WENO_AO_5_3}
  \includegraphics[width=10 cm]{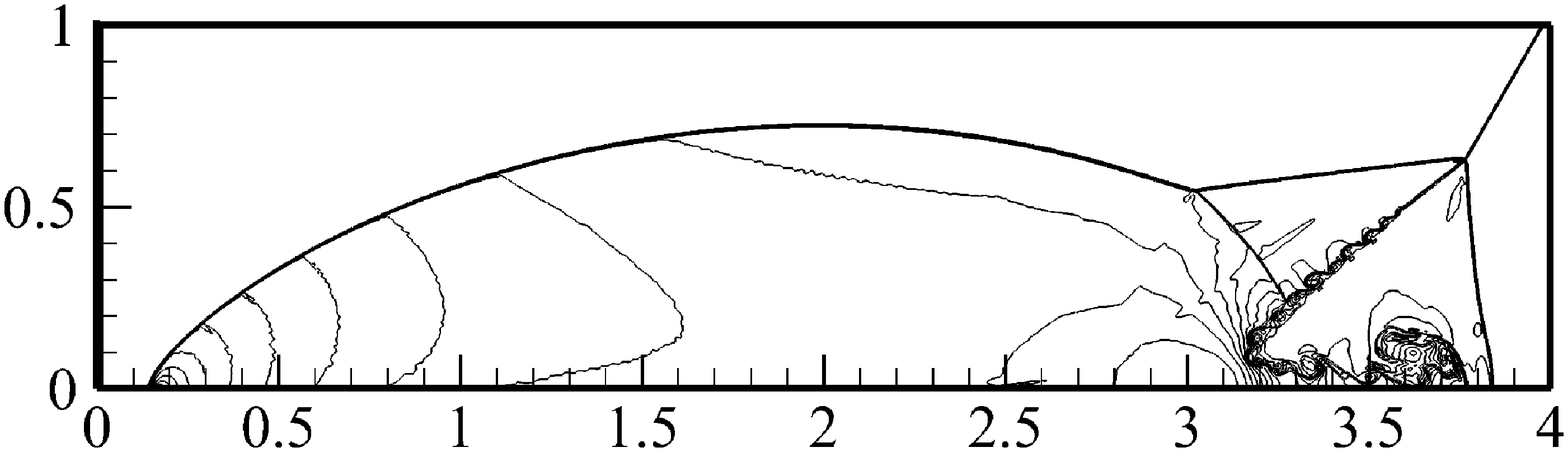}}
  \subfigure[WENO-AOA(5,3)]{
  \label{FIG:DMR_WENO_AOA_5_3}
  \includegraphics[width=10 cm]{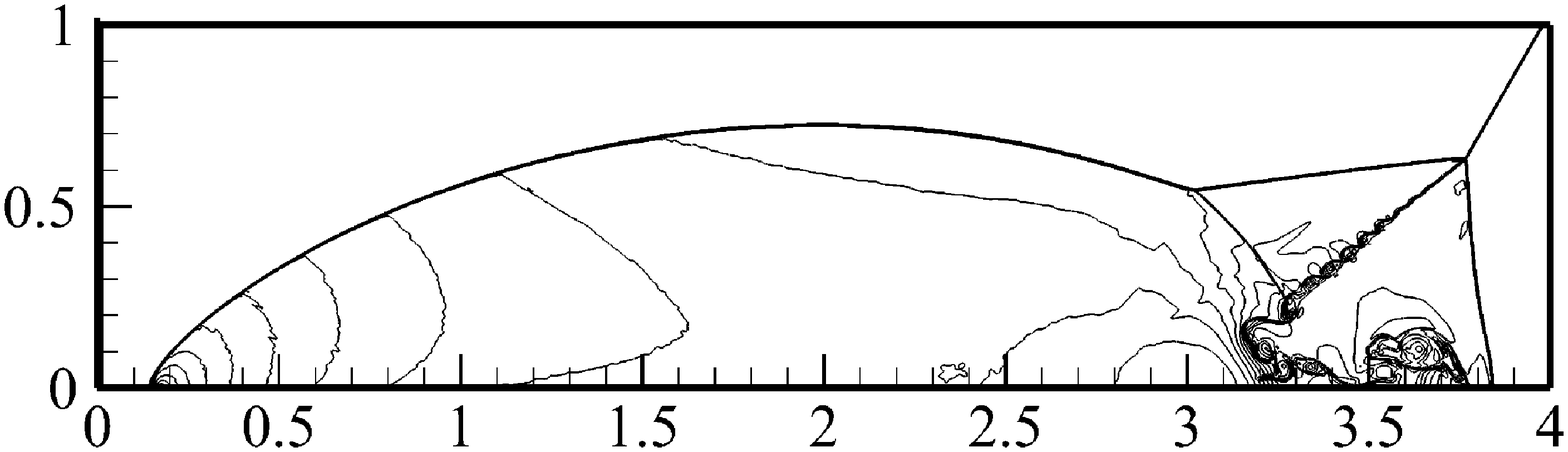}}
  \caption{The entire view of the density contours of double Mach reflection problem at $t=0.28$ calculated by WENO-AO(5,3) and WENO-AOA(5,3) with $1601\times401$ mesh points.}
\label{FIG:DMR}
\end{figure}

\begin{figure}
  \centering
  \subfigure[WENO-AO(5,3)]{
  \label{FIG:DMR_WENO_AO_5_3}
  \includegraphics[width=5.5 cm]{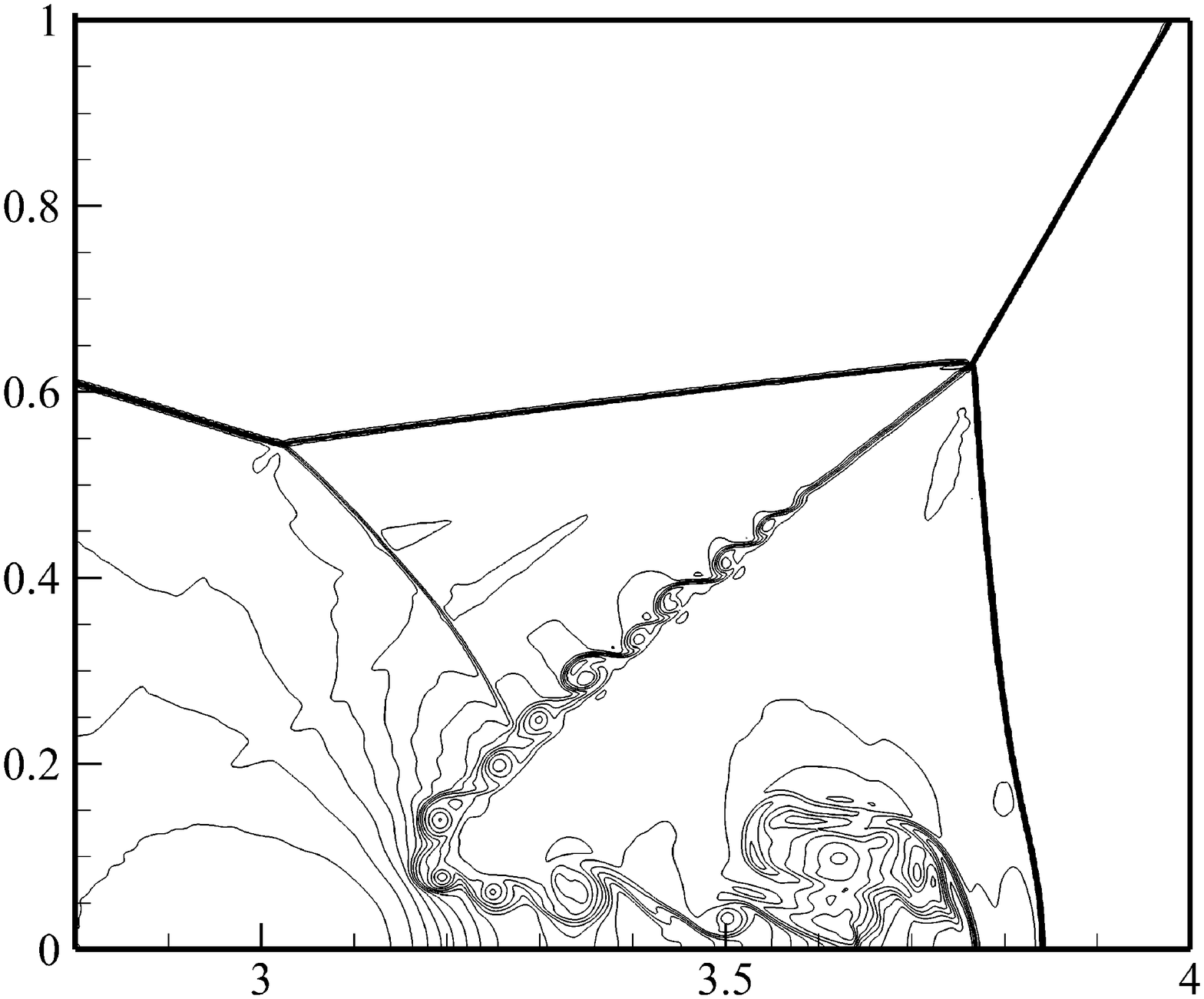}}
  \subfigure[WENO-AOA(5,3)]{
  \label{FIG:DMR_WENO_AOA_5_3}
  \includegraphics[width=5.5 cm]{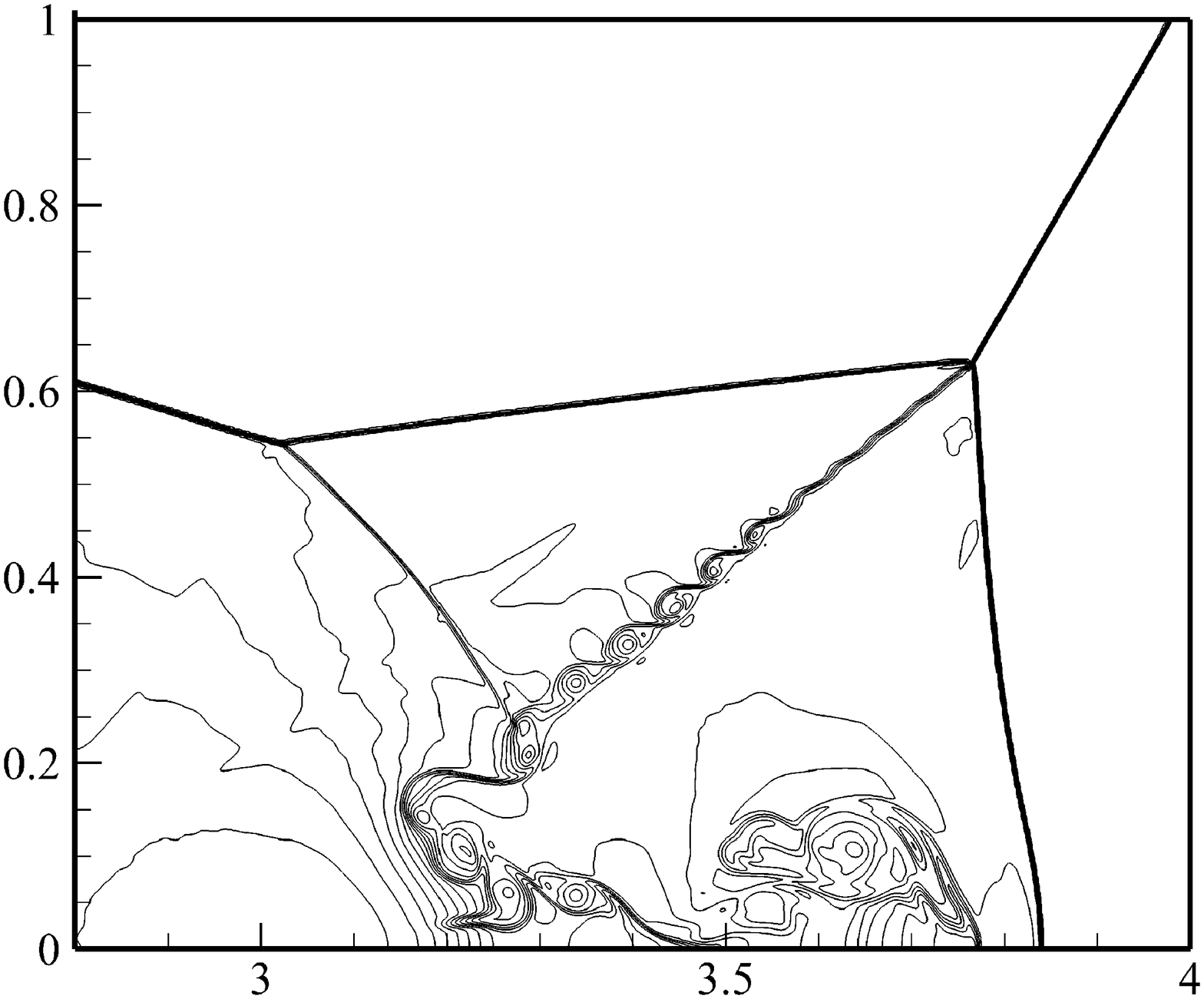}}
  \caption{The enlarged view of the density contours of double Mach reflection problem at $t=0.28$
  calculated by WENO-AO(5,3) and WENO-AOA(5,3) with $1601\times401$ mesh points.}
\label{FIG:DMR_enlarge}
\end{figure}

This problem is governed by the two-dimensional compressible Euler equations
equipped with certain initial conditions and boundary conditions.
It was originally proposed by Woodward and Colella \cite{Woodward1984JCP}
and was widely used to test the performance of numerical schemes
for sophisticated structures and strong shocks.
The specific settings can be found in Woodward and Colella \cite{Woodward1984JCP} and many other papers,
so we omit them for saving space in this note.
We run this problem to $t=0.28$.
Fig. \ref{FIG:DMR} shows the entire view of the density contours calculated
by WENO-AO(5,3) and WENO-AOA(5,3) with $1601\times401$ mesh points.
The overall profiles are the same,
but the details of the reflecting zone are distinct as shown by Fig. \ref{FIG:DMR_enlarge}.
The difference is attributed to the intrinsical chaos characteristic of Kelvin{Helmholtz (KH) instabilities.

\section{Conclusions}
We proposed an alternative WENO-AO reconstruction which has a more concise form than the original WENO-AO reconstruction.
In the benchmark tests, WENO-AOA performs similarly with WENO-AO.
The main advantage of the WENO-AOA is the strict convexity
that may make the reconstruction relatively reliable in general cases.

\section*{Acknowledgement}
H. S. would like to acknowledge the financial support of the National Natural Science Foundation of China (Contract No. 11901602).

\bibliography{mybibfile}

\end{document}